\documentclass[12pt]{article}
\usepackage{amsmath,amssymb,amsthm}
\usepackage[margin=3.2cm]{geometry}
\usepackage{titlesec,hyperref}
\usepackage{fancyhdr}  

\newtheorem{theo}{Theorem}[section]
\newtheorem{lemm}[theo]{Lemma}
\newtheorem{defi}[theo]{Definition}

\newtheorem{rema}[theo]{Remark}
\numberwithin{equation}{section}

\begin{document}
\title{Local existence and uniqueness for the non-resistive MHD equations in homogeneous Besov spaces }

\author{ \footnote{Email Address: lijl29@mail2.sysu.edu.cn (J. Li),\quad  tanwenkeybfq@163.com (W. Tan), \quad  mcsyzy@mail.sysu.edu.cn (Z. Yin, Corresponding author).}
Jinlu $\mbox{Li}^{1,2}$   \quad Wenke $\mbox{Tan}^3$  \quad and \quad Zhaoyang $\mbox{Yin}^{1,4}$ \\
 \small $^1\mbox{Department}$ of Mathematics, Sun Yat-sen University, Guangzhou, 510275, China\\
   \small $^2\mbox{School}$  of Mathematics and Computer Sciences, \\
\small Gannan Normal University, Ganzhou 341000, China\\
   \small $^3\mbox{Department}$ of Mathematics and Computer Science, Hunan Normal University,\\
  \small Changsha 410006, China\\
\small $^4\mbox{Faculty}$ of Information Technology,
\\ \small Macau University of Science and Technology, Macau, China}

\date{}
\maketitle

\begin{abstract}
In the paper, we consider the Cauchy problem of the non-resistive MHD equations in homogeneous Besov spaces. We prove the local existence and uniqueness of the solution to the non-resistive MHD equations by using the iterative scheme and compactness arguments. Our obtained result improves considerably the recent results in \cite{C,W}.

\vspace*{1em}


\vspace*{1em}

\end{abstract}

\section{Introduction}

\quad \quad In this paper, we consider the following non-resistive MHD equations $(d\geq2)$:
\begin{equation}\label{1.1}\begin{cases}
\partial_tu+u\cdot\nabla u-\Delta u+\nabla P=B\cdot \nabla B, \quad (t,x)\in \mathbb{R}^+\times \mathbb{R}^d,\\
\partial_tB+u\cdot \nabla B=B\cdot\nabla u,\quad (t,x)\in \mathbb{R}^+\times \mathbb{R}^d,\\
\mathrm{div} u=\mathrm{div} B=0, \quad (t,x)\in \mathbb{R}^+\times \mathbb{R}^d,\\
(u,B)|_{t=0}=(u_0,B_0), \quad x\in \mathbb{R}^d,
\end{cases}\end{equation}
where the unknowns are the vector fields $u=(u^1,u^2,\cdots,u^d)$, $B=(B^1,B^2,\cdots,B^d)$ and the scalar function $P$. Here, $u$ and $B$ are the velocity and magnetic, respectively, while $P$ denotes the pressure.

Fefferman et al. showed local-in-time existence of strong solutions to \eqref{1.1} in $\mathbb{R}^d,\ d=2,3$ with the initial data $(u_0,B_0)\in H^s(\mathbb{R}^d)\times H^s(\mathbb{R}^d),\ s>\frac d2$ in \cite{F} and $(u_0,B_0)\in H^{s-1-\varepsilon}(\mathbb{R}^d)\times H^s(\mathbb{R}^d),\ s>\frac d2,\ 0<\varepsilon<1$ in \cite{F1}. Chemin et al. in \cite{C} proved the local existence of solutions to \eqref{1.1} in $\mathbb{R}^d,\ d=2,3$ with the initial data $(u_0,B_0)\in \mathfrak{B}^{\frac d2-1}_{2,1}(\mathbb{R}^d)\times \mathfrak{B}^{\frac d2}_{2,1}(\mathbb{R}^d)$ and also proved the corresponding solution is unique in 3D case.  Recently, Wan in \cite{W} obtained the uniqueness of the solution in the 2D case by using mixed space-time Besov spaces.

 However, the local existence and uniqueness of the solution for the Cauchy problem of the non-resistive MHD equations in homogeneous Besov spaces has not been studied yet. Whether or not the solution for the Cauchy problem of the non-resistive MHD equations exists locally in time and is unique in homogeneous Besov spaces  is an open problem which was proposed by Chemin et al. in \cite{C}. In the paper, our aim is to solve this open problem by establishing the local existence and uniqueness of the solution for the Cauchy problem \eqref{1.1} in homogeneous Besov spaces.

 For convenience, we transform the system \eqref{1.1} into an equivalent form of compressible type. By using $\mathrm{div} u=\mathrm{div} B=0$, we have
\begin{align*}
u\cdot\nabla u=\mathrm{div}(u\otimes u), \quad B\cdot \nabla B=\mathrm{div}(B\otimes B), \quad B\cdot\nabla u=\mathrm{div}(u\otimes B).
\end{align*}
Therefore, the system \eqref{1.1} is formally equivalent to the following equations
\begin{equation}\label{1.2}\begin{cases}
\partial_tu-\Delta u=\mathbb{P}\mathrm{div}(B\otimes B-u\otimes u),\\
\partial_tB+u\cdot \nabla B=\mathrm{div}(u\otimes B),\\
(u,B)|_{t=0}=(u_0,B_0),
\end{cases}\end{equation}
where $\mathbb{P}=I+\nabla (-\Delta)^{-1} \mathrm{div}$ and $\mathrm{div}u_0=\mathrm{div}B_0=0$.

To solve \eqref{1.2}, the main difficulty is that the system is only partially parabolic, owing to the magnetic equation which is of hyperbolic type. This precludes any attempt to use the Banach fixed point theorem in a suitable space. Therefore, we would like to present a general functional framework to deal with the local existence and uniqueness of the solution for the compressible fluids of \eqref{1.2} in the homogeneous Besov spaces.

Our main result can be stated as follows:
\begin{theo}\label{th1.1}
Let $d\geq 2$. Assume that the initial data $u_0\in\dot{\mathfrak{B}}^{\frac dp-1}_{p,1}(\mathbb{R}^d)$, and $B_0\in\dot{\mathfrak{B}}^{\frac dp}_{p,1}(\mathbb{R}^d)$. Then there exists a positive time $T$ such that \\

(a) Local existence: If $p\in[1,2d]$, then the system \eqref{1.1} has a solution $(u,B)\in E^p_T$ with
\begin{align*}
E^p_T\triangleq \Big(\mathcal{C}([0,T];\dot{\mathfrak{B}}^{\frac dp-1}_{p,1}(\mathbb{R}^d))\cap L^1_T(\dot{\mathfrak{B}}^{\frac dp+1}_{p,1}(\mathbb{R}^d))\Big)^d\times \Big(\mathcal{C}([0,T];\dot{\mathfrak{B}}^{\frac dp}_{p,1}(\mathbb{R}^d))\Big)^d.
\end{align*}

(b) Uniqueness: If $p\in[1,2d]$, then the uniqueness of the solution holds in $E^p_T$.
\end{theo}

\begin{rema}
Note that when $d=3$, one has $\mathfrak{B}^{\frac 12}_{2,1}(\mathbb{R}^3)\subset\dot{\mathfrak{B}}^{\frac 12}_{2,1}(\mathbb{R}^3)$ and $\mathfrak{B}^{\frac 32}_{2,1}(\mathbb{R}^3)\subset\dot{\mathfrak{B}}^{\frac 32}_{2,1}(\mathbb{R}^3)$. This shows that our obtained result in Theorem 1.1 improves considerably the corresponding result in \cite{C}.
\end{rema}

\begin{rema}
For $d=2$, we also obtain a new result compared with the recent result in \cite{C}. Since $\mathfrak{B}^{0}_{2,1}(\mathbb{R}^2)\subset \dot{\mathfrak{B}}^{0}_{2,1}(\mathbb{R}^2)$ and $\dot{\mathfrak{B}}^{1}_{2,1}(\mathbb{R}^2)\subset \mathfrak{B}^{1}_{2,1}(\mathbb{R}^2)$, then our obtained result in Theorem 1.1 and the corresponding result in \cite{C} don't contain each other.
\end{rema}

\begin{rema}
When $d=2$, scrutinizing our proof of the uniqueness (in Theorem 1.1) for \eqref{1.1} in homogeneous Besov spaces, we can also gain the uniqueness for \eqref{1.1} in nonhomogeneous Besov spaces. This implies that our obtained uniqueness result in Theorem 1.1 covers the recent result in \cite{W}, and our proof is more brief than that in \cite{W}.
\end{rema}

The paper is organized as follows. In Section 2, we recall the Littlewood-Paley theory and give some properties of homogeneous Besov spaces. In Section 3, we prove the local existence and the uniqueness of the solution to the system \eqref{1.2}.
\bigskip
\\
\noindent\textbf{Notations.} In the following, we denote by $\langle\cdot,\cdot\rangle$ the action between $\mathcal{S'}(\mathbb{R}^d)$ and $\mathcal{S}(\mathbb{R}^d)$. Given a Banach space $X$, we denote its norm by $\|\cdot\|_{X}$. Since all spaces of functions are over $\mathbb{R}^d$, for simplicity, we drop  $\mathbb{R}^d$ in our notations of function spaces if there is no ambiguity.

\section{Littlewood-Paley Analysis}

\quad \quad In this section, we first recall the Littlewood-Paley theory, the definition of homogeneous Besov spaces and some useful properties. Then, we state some applications in the linear transport equation and the heat conductive equation.

First, let us introduce the Littlewood-Paley decomposition. Choose a radial function $\varphi\in \mathcal{S}(\mathbb{R}^d)$ supported in $\tilde{\mathcal{C}}=\{\xi\in\mathbb{R}^d,\frac34\leq \xi\leq \frac83\}$ such that
\begin{align*}
\sum_{j\in \mathbb{Z}}\varphi(2^{-j}\xi)=1 \quad \mathrm{for} \ \mathrm{all} \ \xi\neq0.
\end{align*}
The frequency localization operator $\dot{\Delta}_j$ and $\dot{S}_j$ are defined by
\begin{align*}
\dot{\Delta}_jf=\varphi(2^{-j}D)f=\mathcal{F}^{-1}(\varphi(2^{-j}\cdot)\mathcal{F}f), \quad \dot{S}_jf=\sum_{k\leq j-1}\dot{\Delta}_kf \quad \mathrm{for} \quad j\in\mathbb{Z}.
\end{align*}
With a suitable choice of $\varphi$, one can easily verify that
\begin{align*}
\dot{\Delta}_j\dot{\Delta}_kf=0 \quad \mathrm{if} \quad |j-k|\geq2, \quad \dot{\Delta}_j(\dot{S}_{k-1}f\dot{\Delta}_kf)=0 \quad  \mathrm{if} \quad  |j-k|\geq5.
\end{align*}
Next we recall Bony's decomposition from \cite{B.C.D}:
\begin{align*}
uv=\dot{T}_uv+\dot{T}_vu+\dot{R}(u,v),
\end{align*}
with
\begin{align*}
\dot{T}_uv=\sum_{j\in\mathbb{Z}}\dot{S}_{j-1}\dot{\Delta}_jv, \quad \quad \dot{R}(u,v)=\sum_{j\in\mathbb{Z}}\dot{\Delta}_ju\widetilde{\Delta}_jv, \quad \quad \widetilde{\Delta}_jv=\sum_{|j'-j|\leq 1}\dot{\Delta}_{j'}v.
\end{align*}

The following Bernstein lemma will be stated as follows (see \cite{B.C.D}):
\begin{lemm}
Let $1\leq p\leq q\leq \infty$  and $\mathcal{B}$ be a ball and $\mathcal{C}$ a ring of $\mathbb{R}^d$. Assume that $f\in L^p$, then for any $\alpha\in\mathbb{N}^d$, there exists a constant $C$ independent of $f$, $j$ such that
\begin{align*}
&&\mathrm{Supp} \,\hat{f}\subset\lambda \mathcal{B}\Rightarrow\sup_{|\alpha|=k}\|\partial^{\alpha}f\|_{L^q}\le C^{k+1}\lambda^{k+d(\frac1p-\frac1q)}\|f\|_{L^p},\\
&&\mathrm{Supp} \,\hat{f}\subset\lambda \mathcal{C}\Rightarrow C^{-k-1}\lambda^k\|f\|_{L^p}\le\sup_{|\alpha|=k}\|\partial^{\alpha}f\|_{L^p}
\le C^{k+1}\lambda^{k}\|f\|_{L^p}.
\end{align*}
\end{lemm}

Now, we will introduce the definition of the homogeneous Besov space. We denote the space $\mathcal{Z}'(\mathbb{R}^d)$ by the dual space of $\mathcal{Z}(\mathbb{R}^d)=\{f\in \mathcal{S}(\mathbb{R}^d);D^{\alpha}\hat{f}(0)=0;\ \forall \alpha\in \mathbb{N}^d\}$, which can be identified by the quotient space of $\mathcal{S}'(\mathbb{R}^d)/\mathcal{P}$ with the polynomials space $\mathcal{P}$. The formal equality $f=\sum\limits_{j\in \mathbb{Z}}\dot{\Delta}_jf$ holds true for $f\in\mathcal{Z}'(\mathbb{R}^d)$ and is called the homogenous Littlewood-Paley decomposition.

The operators $\dot{\Delta}_j$ help us recall the definition of the homogenous Besov space (see \cite{B.C.D,T})

\begin{defi}
Let $s\in \mathbb{R}$, $1\leq p,r\leq \infty$. The homogeneous Besov space $\dot{\mathfrak{B}}^s_{p,r}$ is defined by
\begin{align*}
\dot{\mathfrak{B}}^s_{p,r}=\{f\in \mathcal{Z}'(\mathbb{R}^d):||f||_{\dot{\mathfrak{B}}^s_{p,r}}<+\infty\},
\end{align*}
where
\begin{align*}
||f||_{\dot{\mathfrak{B}}^s_{p,r}}\triangleq \Big|\Big|(2^{ks}||\dot{\Delta}_k f||_{L^p})_{k}\Big|\Big|_{\ell^r}.
\end{align*}
\end{defi}
It is easy to check that $||\mathbb{P}f||_{\dot{\mathfrak{B}}^s_{p,r}}\leq C||f||_{\dot{\mathfrak{B}}^s_{p,r}}$ for some positive constant $C$.
We also need to use Chemin-Lerner type Besov spaces introduced in \cite{B.C.D},
\begin{defi}
Let $s\in\mathbb{R}$, $1\leq p,q,r\leq\infty$ and $T\in(0,\infty]$. The functional space $\tilde{L}^q_T(\dot{\mathfrak{B}}^s_{p,r})$ is defined as the set of all the distributions $f(t)$ satisfying
\begin{align*}
||f||_{\tilde{L}^q_T(\dot{\mathfrak{B}}^s_{p,r})}\triangleq \Big|\Big|\big(2^{ks}||\dot{\Delta}_kf(t)||_{L^q_T(L^p)}\big)_k \Big|\Big|_{\ell^r}<+\infty.
\end{align*}
\end{defi}
By Minkowski's inequality, it is easy to find that
\begin{align*}
||f||_{\tilde{L}^q_T(\dot{\mathfrak{B}}^s_{p,r})}\leq ||f||_{L^q_T(\dot{\mathfrak{B}}^s_{p,r})} \quad \mathrm{if} \quad q\leq r, \quad \quad ||f||_{\tilde{L}^q_T(\dot{\mathfrak{B}}^s_{p,r})}\geq ||f||_{L^q_T(\dot{\mathfrak{B}}^s_{p,r})} \quad \mathrm{if} \quad q\geq r.
\end{align*}

We state the following logarithmic interpolation inequality which will be useful in the sequel (see \cite{D1}).
\begin{lemm}(see Proposition 2.8 in \cite{D1})\label{le2.3}
Let $s\in\mathbb{R}$. Then, for any $1\leq p,q\leq \infty$ and $0<\varepsilon\leq 1$, we have
\begin{align*}
||f||_{\tilde{L}^q_T(\dot{\mathfrak{B}}^s_{p,1})}\leq C\frac{||f||_{\tilde{L}^q_T(\dot{\mathfrak{B}}^s_{p,\infty})}}{\varepsilon}\log(e+\frac{||f||_{\tilde{L}^q_T(\dot{\mathfrak{B}}^{s-\varepsilon}_{p,\infty})}+
||f||_{\tilde{L}^q_T(\dot{\mathfrak{B}}^{s+\varepsilon}_{p,\infty})}}{||f||_{\tilde{L}^q_T(\dot{\mathfrak{B}}^s_{p,\infty})}}).
\end{align*}
\end{lemm}

Next, we give the important product acts on homogenous Besov spaces or Chemin-Lerner type Besov spaces by collecting some useful lemmas from \cite{C.M.Z,D}.
\begin{lemm}(see Proposition 1.8 in \cite{D})\label{le2.4}
Let $1\leq p,q,q_1,q_2\leq \infty$ with $\frac{1}{q_1}+\frac{1}{q_2}=\frac1q$. Then there hold \\
(a) if $s_2\leq \frac dp$, we have
\begin{align*}
||\dot{T}_gf||_{\dot{\mathfrak{B}}^{s_1+s_2-\frac dp}_{p,1}}&\leq C||f||_{\dot{\mathfrak{B}}^{s_1}_{p,1}}||g||_{\dot{\mathfrak{B}}^{s_2}_{p,1}};\\
||\dot{T}_gf||_{\tilde{L}^q_T(\dot{\mathfrak{B}}^{s_1+s_2-\frac dp}_{p,1})}&\leq C||f||_{\tilde{L}^{q_1}_T(\dot{\mathfrak{B}}^{s_1}_{p,1})}||g||_{\tilde{L}^{q_2}_T(\dot{\mathfrak{B}}^{s_2}_{p,1})};
\end{align*}
(b) if $s_1+s_2>d\max\{0,\frac2p-1\}$, we have
\begin{align*}
||\dot{R}(f,g)||_{\dot{\mathfrak{B}}^{s_1+s_2-\frac dp}_{p,1}}&\leq C||f||_{\dot{\mathfrak{B}}^{s_1}_{p,1}}||g||_{\dot{\mathfrak{B}}^{s_2}_{p,1}};\\
||\dot{R}(f,g)||_{\tilde{L}^q_T(\dot{\mathfrak{B}}^{s_1+s_2-\frac dp}_{p,1})}&\leq C||f||_{\tilde{L}^{q_1}_T(\dot{\mathfrak{B}}^{s_1}_{p,1})}||g||_{\tilde{L}^{q_2}_T(\dot{\mathfrak{B}}^{s_2}_{p,1})}.
\end{align*}
\end{lemm}

\begin{lemm}(see Lemma 2.6 in \cite{C.M.Z})\label{le2.5}
Let $s_1,s_2\leq \frac dp$, $s_1+s_2>d\max\{0,\frac2p-1\}$ and $1\leq p,q,q_1,q_2\leq \infty$ with $\frac{1}{q_1}+\frac{1}{q_2}=\frac1q$. Then there hold
\begin{align*}
||fg||_{\dot{\mathfrak{B}}^{s_1+s_2-\frac dp}_{p,1}}&\leq C||f||_{\dot{\mathfrak{B}}^{s_1}_{p,1}}||g||_{\dot{\mathfrak{B}}^{s_2}_{p,1}},\\
||fg||_{\tilde{L}^q_T(\dot{\mathfrak{B}}^{s_1+s_2-\frac dp}_{p,1})}&\leq C||f||_{\tilde{L}^{q_1}_T(\dot{\mathfrak{B}}^{s_1}_{p,1})}||g||_{\tilde{L}^{q_2}_T(\dot{\mathfrak{B}}^{s_2}_{p,1})}.
\end{align*}
\end{lemm}

\begin{lemm}(see Lemma 2.7 in \cite{C.M.Z})\label{le2.6}
Let $s_1\leq \frac dp,\ s_2<\frac dp$, $s_1+s_2\geq d\max\{0,\frac2p-1\}$ and $1\leq p,q,q_1,q_2\leq \infty$ with $\frac{1}{q_1}+\frac{1}{q_2}=\frac1q$. Then there hold
\begin{align*}
||fg||_{\dot{\mathfrak{B}}^{s_1+s_2-\frac dp}_{p,\infty}}&\leq C||f||_{\dot{\mathfrak{B}}^{s_1}_{p,1}}||g||_{\dot{\mathfrak{B}}^{s_2}_{p,\infty}},\\
||fg||_{\tilde{L}^q_T(\dot{\mathfrak{B}}^{s_1+s_2-\frac dp}_{p,\infty})}&\leq C||f||_{\tilde{L}^{q_1}_T(\dot{\mathfrak{B}}^{s_1}_{p,1})}||g||_{\tilde{L}^{q_2}_T(\dot{\mathfrak{B}}^{s_2}_{p,\infty})}.
\end{align*}
\end{lemm}

Finally, we will present the priori estimates of the linear transport equation
\begin{align}\label{3.1}
\partial_tf+v\cdot\nabla f=g,\quad f(0,x)=f_0,
\end{align}
and the heat conductive equation
\begin{align}\label{3.2}
\partial_tu-\Delta u=G, \quad u(0,x)=u_0,
\end{align}
in homogenous Besov spaces. The following estimates will be frequently used in the sequel (see \cite{D2}).
\begin{lemm}(see Proposition 1.8 in \cite{D2})\label{le3.1}
Let $s\in(-d\min\{\frac1p,1-\frac{1}{p}\}-1,1+\frac dp)$, $1\leq p,r\leq \infty$ and $s=1+\frac dp$ if $r=1$. Let $v$ be a vector field such that $\nabla v\in L^1_T(\dot{\mathfrak{B}}^{\frac dp}_{p,r}\cap L^\infty)$ and $\mathrm{div} v=0$. Assume that $f_0\in\dot{\mathfrak{B}}^s_{p,r}$, $g\in L^1_T(\dot{\mathfrak{B}}^s_{p,r})$ and $f$ is the solution of \eqref{3.1}. Then there holds for $t\in[0,T]$,
\begin{align*}
||f||_{\tilde{L}^\infty_t(\dot{\mathfrak{B}}^s_{p,r})}\leq e^{CV(t)}(||f_0||_{\dot{\mathfrak{B}}^s_{p,r}}+\int^t_0e^{-CV(\tau)}||g(\tau)||_{\dot{\mathfrak{B}}^s_{p,r}}\mathrm{d}\tau),
\end{align*}
or
\begin{align*}
||f||_{\tilde{L}^\infty_t(\dot{\mathfrak{B}}^s_{p,r})}\leq e^{CV(t)}(||f_0||_{\dot{\mathfrak{B}}^s_{p,r}}+||g||_{\tilde{L}^1_t(\dot{\mathfrak{B}}^s_{p,r})}),
\end{align*}
where $V(t)=\int^t_0||\nabla v||_{\dot{\mathfrak{B}}^{\frac dp}_{p,r}\cap L^\infty}\mathrm{d} \tau$.
\end{lemm}

\begin{lemm}(see Proposition 1.9 in \cite{D2})\label{le3.2}
Let $s\in \mathbb{R}$ and $1\leq q,q_1,p,r\leq \infty$ with $q_1\leq q$. Assume that $u_0\in \dot{\mathfrak{B}}^s_{p,r}$ and $G\in \tilde{L}^{q_1}_T(\dot{\mathfrak{B}}^{s-2+\frac{2}{q_1}}_{p,r})$. Then \eqref{3.2} has a unique solution $u\in \tilde{L}^{q}_T(\dot{\mathfrak{B}}^{s+\frac2q}_{p,r})$ satisfying
\begin{align*}
||u||_{\tilde{L}^{q}_T(\dot{\mathfrak{B}}^{s+\frac2q}_{p,r})}\leq C(||u_0||_{\dot{\mathfrak{B}}^s_{p,r}}+||G||_{\tilde{L}^{q_1}_T(\dot{\mathfrak{B}}^{s-2+\frac{2}{q_1}}_{p,r})}).
\end{align*}
\end{lemm}

\section{Local existence and uniqueness in homogeneous Besov spaces}
\par
\quad \quad In this section, we will show the local existence and uniqueness of the solution to the system \eqref{1.2} with the initial data in homogeneous Besov spaces. However, because the whole system is not fully parabolic, the strong convergence of the sequence is shown for a weaker norm corresponding to a loss of one derivative. For that reason, the uniqueness issue is tractable by taking advantage of a logarithmic interpolation inequality together with Osgood's lemma.

In the following, we divide the proof of Theorem 1.1 into four steps to prove the local existence and the uniqueness of the solution to the system \eqref{1.2}.\\

\noindent \textbf{Step 1: An iterative scheme. }We will use an iterative scheme to obtain the approximating sequence to the system \eqref{1.2} by combining the linear transport equation and the heat conductive equation. Set $u^n_0\triangleq\dot{S}_nu_0$ and $B^n_0=\dot{S}_nB_0$, and define the first term $(u^0, B^0)$ of the approximating sequence to be
\begin{align*}
u^0\triangleq e^{t\Delta}u^0_0, \quad  B^0\triangleq e^{t\Delta}B^0_0.
\end{align*}
Starting from the above term $(u^0, B^0)$, we define by induction a sequence ($u^{n}, B^{n})_{n\in\mathbb{N}}$ of smooth functions by solving the following linear transport and heat conductive equations:
\begin{equation}\label{4.1}\begin{cases}
\partial_tu^{n+1}-\Delta u^{n+1}=\mathbb{P}\mathrm{div}(-u^n\otimes u^n+B^n\otimes B^n),\\
\partial_tB^{n+1}-u^n\cdot\nabla B^{n+1}=\mathrm{div}(u^n\otimes B^n),\\
(u^{n+1}_0,B^{n+1}_0)=\dot{S}_{n+1}(u_0,B_0).
\end{cases}\end{equation}

\noindent \textbf{Step 2: Uniform estimates. }Taking advantage of Lemmas \ref{le3.1} - \ref{le3.2}, we shall bound the approximating sequence in the expected solution space on some fixed time interval.  That is, for all $T>0$ and $n\in \mathbb{N}$, we have $(u^n,B^n)\in E^p_T$. Now, we claim that there exists some $T > 0$ independent of $n$ such that the solution $(u^n,B^n)_{n\in \mathbb{N}}$ satisfies the following inequalities for some positive constants $C_0>1$ and $\eta<1$ (to be determined later):
\begin{align*}
&(H_1):\quad ||u^n||_{\tilde{L}^\infty_T(\dot{\mathfrak{B}}^{\frac dp-1}_{p,1})}+||B^n||_{\tilde{L}^\infty_T(\dot{\mathfrak{B}}^{\frac dp}_{p,1})}
\leq C_0E_0,\\
&(H_2):\quad ||u^n||_{\tilde{L}^1_T(\dot{\mathfrak{B}}^{\frac dp+1}_{p,1})}+||u^n||_{\tilde{L}^2_T(\dot{\mathfrak{B}}^{\frac dp}_{p,1})}\leq \eta,
\end{align*}
where $E_0=||u_0||_{\dot{\mathfrak{B}}^{\frac dp-1}_{p,1}}+||B_0||_{\dot{\mathfrak{B}}^{\frac dp}_{p,1}}$. Now, we suppose that $T$ satisfies the following inequality:
\begin{align}\label{4}
||e^{t\Delta}u_0||_{\tilde{L}^1_T(\dot{\mathfrak{B}}^{\frac dp+1}_{p,1})}+||e^{t\Delta}u_0||_{\tilde{L}^2_T(\dot{\mathfrak{B}}^{\frac dp}_{p,1})}\leq \eta^2.
\end{align}
It is easy to check the conditions $(H_1)$ - $(H_2)$ hold true for $n=0$ by \eqref{4}. In what follows, we will show that if the conditions $(H_1)$ - $(H_2)$ hold true for $n$, then they hold true for $n+1$. First of all, we get by Lemma \ref{le2.5},
\begin{align}\label{44}
&\quad \ ||\mathrm{div}(-u^n\otimes u^n+B^n\otimes B^n)||_{L^1_T(\dot{\mathfrak{B}}^{\frac dp-1}_{p,1})}\leq C||-u^n\otimes u^n+B^n\otimes B^n)||_{L^1_T(\dot{\mathfrak{B}}^{\frac dp}_{p,1})} \nonumber
\\&\leq C||u^n||_{\tilde{L}^2_T(\dot{\mathfrak{B}}^{\frac dp}_{p,1})}||u^n||_{\tilde{L}^2_T(\dot{\mathfrak{B}}^{\frac dp}_{p,1})}+C||B^n||_{\tilde{L}^\infty_T(\dot{\mathfrak{B}}^{\frac dp}_{p,1})}||B^n||_{\tilde{L}^1_T(\dot{\mathfrak{B}}^{\frac dp}_{p,1})}.
\end{align}
It follows from Lemma \ref{le3.2} and \eqref{44} that
\begin{align}\label{4.2}
&\quad \ ||u^{n+1}||_{\tilde{L}^\infty_T(\dot{\mathfrak{B}}^{\frac dp-1}_{p,1})} \nonumber
\\&\leq C||u^{n+1}_0||_{\dot{\mathfrak{B}}^{\frac dp-1}_{p,1}}+C||\mathrm{div}(-u^n\otimes u^n+B^n\otimes B^n)||_{L^1_T(\dot{\mathfrak{B}}^{\frac dp-1}_{p,1})}\nonumber
\\&\leq C||u_0||_{\dot{\mathfrak{B}}^{\frac dp-1}_{p,1}}+C||u^n||^2_{\tilde{L}^2_T(\dot{\mathfrak{B}}^{\frac dp}_{p,1})}+CT||B^n||^2_{\tilde{L}^\infty_T(\dot{\mathfrak{B}}^{\frac dp}_{p,1})}.
\end{align}
By Lemma \ref{le2.5} and Lemma \ref{le3.1}, we have
\begin{align}\label{4.3}
||B^{n+1}||_{\tilde{L}^\infty_T(\dot{\mathfrak{B}}^{\frac dp}_{p,1})}&\leq Ce^{CU^n(T)}(||B^n_0||_{\dot{\mathfrak{B}}^{\frac dp}_{p,1}}+||B^n\cdot \nabla u^n||_{\tilde{L}^1_T(\dot{\mathfrak{B}}^{\frac dp}_{p,1})})\nonumber
\\& \leq Ce^{CU^n(T)}(||B_0||_{\dot{\mathfrak{B}}^{\frac dp}_{p,1}}+||B^n||_{\tilde{L}^\infty_T(\dot{\mathfrak{B}}^{\frac dp}_{p,1})}||u^n||_{\tilde{L}^1_T(\dot{\mathfrak{B}}^{\frac dp+1}_{p,1})}),
\end{align}
where $U^n(t)=\int^t_0||u^n||_{\dot{\mathfrak{B}}^{\frac dp+1}_{p,1}}\mathrm{d} \tau$. Then, combining \eqref{4.2}-\eqref{4.3} and the conditions $(H_1)$ - $(H_2)$, we obtain
\begin{align}\label{40}
&\quad \ ||u^{n+1}||_{\tilde{L}^\infty_T(\dot{\mathfrak{B}}^{\frac dp-1}_{p,1})}+||B^{n+1}||_{\tilde{L}^\infty_T(\dot{\mathfrak{B}}^{\frac dp}_{p,1})}\leq Ce^{C\eta}E_0+CC^2_0E^2_0T+C\eta e^{C\eta}C_0E_0+C\eta^2.
\end{align}
Next, we also get by Lemma \ref{le3.2} and \eqref{4} that
\begin{align}\label{41}
&\quad \ ||u^{n+1}||_{\tilde{L}^1_T(\dot{\mathfrak{B}}^{\frac dp+1}_{p,1})}+||u^{n+1}||_{\tilde{L}^2_T(\dot{\mathfrak{B}}^{\frac dp}_{p,1})} \nonumber
\\&\leq C(||e^{t\Delta}u^{n+1}_0||_{\tilde{L}^1_T(\dot{\mathfrak{B}}^{\frac dp+1}_{p,1})}+||e^{t\Delta}u^{n+1}_0||_{\tilde{L}^2_T(\dot{\mathfrak{B}}^{\frac dp}_{p,1})}) +C||\mathrm{div}(-u^n\otimes u^n+B^n\otimes B^n)||_{L^1_T(\dot{\mathfrak{B}}^{\frac dp-1}_{p,1})}
\nonumber \\&\leq C(||e^{t\Delta}u_0||_{\tilde{L}^1_T(\dot{\mathfrak{B}}^{\frac dp+1}_{p,1})}+||e^{t\Delta}u_0||_{\tilde{L}^2_T(\dot{\mathfrak{B}}^{\frac dp}_{p,1})}) + C||u^n||^2_{\tilde{L}^2_T(\dot{\mathfrak{B}}^{\frac dp}_{p,1})}+CT||B^n||^2_{\tilde{L}^\infty_T(\dot{\mathfrak{B}}^{\frac dp}_{p,1})}\nonumber
\\&\leq C\eta^2+CC^2_0E^2_0T.
\end{align}
Therefore, if choosing $\eta<\min\{\frac{1}{4C},E_0\}$, $C_0>8(C+1)$ and $T<\frac{\eta^2}{C^2_0E^2_0+1}$, then we infer from \eqref{40}-\eqref{41} that
\begin{align*}
&||u^{n+1}||_{\tilde{L}^\infty_T(\dot{\mathfrak{B}}^{\frac dp-1}_{p,1})}+||B^{n+1}||_{\tilde{L}^\infty_T(\dot{\mathfrak{B}}^{\frac dp}_{p,1})}\leq C_0E_0,\\
& ||u^{n+1}||_{\tilde{L}^1_T(\dot{\mathfrak{B}}^{\frac dp+1}_{p,1})}+||u^{n+1}||_{\tilde{L}^2_T(\dot{\mathfrak{B}}^{\frac dp}_{p,1})}\leq \eta.
\end{align*}
This implies the conditions $(H_1)$ - $(H_2)$ hold true for $n+1$. Therefore, we deduce that the approximate sequence $(u^n,B^n)_{n\in \mathbb{N}}$ is uniformly bounded in $E^p_T$ independent of $n$.\\

\noindent \textbf{Step 3: Existence of a solution. }Now, we will use the compactness argument in Besov spaces for the approximating sequence $(u^n,B^n)_{n\in \mathbb{N}}$ to get some solution $(u,B)$ which also satisfies the system \eqref{1.2} in the sense of distributions. Since $u^n$ is uniformly bounded in $L^\infty_T(\dot{\mathfrak{B}}^{\frac dp-1}_{p,1})\cap L^1_T(\dot{\mathfrak{B}}^{\frac dp+1}_{p,1})$, the interpolation inequality yields that $u^n$ is uniformly bounded in $L^q_T(\dot{\mathfrak{B}}^{\frac dp-1+\frac{2}{q}}_{p,1})$ for $1\leq q\leq \infty$. Then, we get by Lemma \ref{le2.5} and $\mathrm{div} u^n=0$,
\begin{align*}
&||u^n\cdot\nabla B^{n+1}-\mathrm{div}(u^n\otimes B^n)||_{L^2_T(\dot{\mathfrak{B}}^{\frac dp-1}_{p,1})}=||\mathrm{div}(B^{n+1}\otimes u^n-u^n\otimes B^n)||_{L^2_T(\dot{\mathfrak{B}}^{\frac dp-1}_{p,1})}\\&\leq C||u^n||_{L^2_T(\dot{\mathfrak{B}}^{\frac dp}_{p,1})}||B^{n+1}||_{L^\infty_T(\dot{\mathfrak{B}}^{\frac dp}_{p,1})}+C||u^n||_{L^2_T(\dot{\mathfrak{B}}^{\frac dp}_{p,1})}||B^{n}||_{L^\infty_T(\dot{\mathfrak{B}}^{\frac dp}_{p,1})}.
\end{align*}
We infer that $\partial_tB^{n+1}$ is uniformly bounded in $L^2_T(\dot{\mathfrak{B}}^{\frac dp-1}_{p,1})$. On the other hand, by Lemma \ref{le2.4}, we have
\begin{align*}
&||\mathrm{div}(u^n\otimes u^n)||_{\dot{\mathfrak{B}}^{\frac dp-\frac32}_{p,1}} \leq C||u^n||_{\dot{\mathfrak{B}}^{\frac dp-1}_{p,1}}||u^n||_{\dot{\mathfrak{B}}^{\frac dp+\frac12}_{p,1}},
\\&||\mathrm{div}(B^n\otimes B^n)||_{\dot{\mathfrak{B}}^{\frac dp-1}_{p,1}} \leq C||B^n||_{\dot{\mathfrak{B}}^{\frac dp}_{p,1}}||B^n||_{\dot{\mathfrak{B}}^{\frac dp}_{p,1}}.
\end{align*}
Then, from the second equation of the system (3.1) and the above inequalities, we infer that $(\partial_tu^n)_{n\in \mathbb{N}}$ is uniformly bounded in $L^\frac43_T(\dot{\mathfrak{B}}^{\frac dp-\frac32}_{p,1}+\dot{\mathfrak{B}}^{\frac dp-1}_{p,1})$.

Let $(\chi_j)_{j\in\mathbb{N}}$ be a sequence of smooth functions with value in $[0,1]$ supported in the ball $B(0,j+1)$ and equal to 1 on $B(0,j)$. The above argument ensures that $(B^n)_{n\in\mathbb{N}}$ is uniformly bounded in $\mathcal{C}^{\frac12}([0,T];\dot{\mathfrak{B}}^{\frac dp-1}_{p,1})\cap \mathcal{C}([0,T];\dot{\mathfrak{B}}^{\frac dp}_{p,1})$, and $(u^n)_{n\in\mathbb{N}}$ is uniformly bounded in $\mathcal{C}^{\frac14}([0,T];\dot{\mathfrak{B}}^{\frac dp-\frac32}_{p,1}+\dot{\mathfrak{B}}^{\frac dp-1}_{p,1})\cap \mathcal{C}([0,T];\dot{\mathfrak{B}}^{\frac dp-1}_{p,1})$. Then, By Proposition 2.93 in \cite{B.C.D} and the embedding relation $\mathfrak{B}^{\frac dp-1}_{p,1}\hookrightarrow \mathfrak{B}^{\frac dp-\frac32}_{p,1}$, we get that for any $j\in\mathbb{N}$, $(\chi_jB^n)_{n\in\mathbb{N}}$ is uniformly bounded in $\mathcal{C}^{\frac12}([0,T];\mathfrak{B}^{\frac dp-1}_{p,1})\cap \mathcal{C}([0,T];\mathfrak{B}^{\frac dp}_{p,1})$, and $(\chi_ju^n)_{n\in\mathbb{N}}$ is uniformly bounded in $\mathcal{C}^{\frac14}([0,T];\mathfrak{B}^{\frac dp-\frac32}_{p,1})\cap \mathcal{C}([0,T];\mathfrak{B}^{\frac dp-1}_{p,1})$. Now, according to Theorem 2.94 in \cite{B.C.D}, the map $z\rightarrow \chi_jz$ is compact from $\mathfrak{B}^{\frac dp}_{p,1}$ to $\mathfrak{B}^{\frac dp-1}_{p,1}$ and $\mathfrak{B}^{\frac dp-1}_{p,1}$ to $\mathfrak{B}^{\frac dp-\frac32}_{p,1}$. Thus, by applying Ascoli's theorem (see Theorem 2.1 in \cite{Feireisl}) and Cantor's diagonal process, there exists some function $(u_j,B_j)$ such that for any $j\in \mathbb{N}$, $\chi_ju$ tends to $u_j$ and $\chi_jB$ tends to $B_j$. As $\chi_j\chi_{j+1}=\chi_{j}$, we have, in addition, $u_j=\chi_ju_{j+1}$ and $B_j=\chi_jB_{j+1}$. From that, we can easily deduce that there exists $(u,B)$ such that for all $\chi\in \mathcal{D}(\mathbb{R}^d)$,
\begin{equation}\label{42}\begin{cases}
 \chi B^n\rightarrow \chi B \quad \mathrm{in} \quad \mathcal{C}([0,T];\mathfrak{B}^{\frac dp-1}_{p,1}),\\
 \chi u^n\rightarrow \chi u \quad \mathrm{in} \quad \mathcal{C}([0,T];\mathfrak{B}^{\frac dp-\frac32}_{p,1}),
\end{cases}\end{equation}
as $n$ tends to $\infty$ (up to a subsequence). By the interpolation, we also have
\begin{equation}\label{43}\begin{cases}
 \chi B^n\rightarrow \chi B \quad \mathrm{in} \quad \ \mathcal{C}([0,T];\mathfrak{B}^{\frac dp-s}_{p,1}), \quad \mathrm{for} \  \mathrm{all} \quad 0< s\leq1,\\
 \chi u^n\rightarrow \chi u \quad \mathrm{in} \quad L^1([0,T];\mathfrak{B}^{\frac dp+s}_{p,1}), \quad \mathrm{for} \  \mathrm{all} \ -\frac32\leq s<1.
\end{cases}\end{equation}
Combining the uniform bounds which we have proved in Step 2 and the Fatou property for Besov spaces, we readily get
\begin{align}\label{99}
(u,B)\in (\tilde{L}^\infty_T(\dot{\mathfrak{B}}^{\frac dp-1}_{p,1})\cap \tilde{L}^2_T(\dot{\mathfrak{B}}^{\frac dp}_{p,1}))^d\times \Big(\tilde{L}^\infty_T(\dot{\mathfrak{B}}^{\frac dp}_{p,1})\Big)^d.
\end{align}
With \eqref{42}-\eqref{43}, it is a routine process to verify that $(u,B)$ satisfies the system \eqref{1.2} in the sense of distributions. For simplicity, we only check the nonlinear term $\mathrm{div}(u^n\otimes B^n)$. For any $\phi\in (\mathcal{D}(\mathbb{R}^d))^d$, there exists some $\chi\in \mathcal{D}(\mathbb{R}^d)$ such that $\phi=\phi\chi$. Therefore, letting $p'$ be the H\"{o}lder conjugate of $p$, i.e. $\frac{1}{p}+\frac{1}{p'}=1$, we infer from Lemma \ref{le2.6} that
\begin{align*}
& \quad \ \big|\int^T_0\langle \mathrm{div}(u^n\otimes B^n-u\otimes B),\phi\rangle\mathrm{d} t\big|=\big|\int^T_0\langle \mathrm{div}\big((\chi u^n)\otimes (\chi B^n)-(\chi u)\otimes (\chi B)\big),\phi\rangle\mathrm{d} t\big|
\\&\leq\big|\int^T_0\langle \mathrm{div}\big((\chi u^n-\chi u)\otimes (\chi B^n)\big),\phi\rangle\mathrm{d} t\big|+
\big|\int^T_0\langle \mathrm{div}\big((\chi u)\otimes (\chi B^n-\chi B)\big),\phi\rangle\mathrm{d} t\big|
\\&\leq C||\chi u^n-\chi u||_{L^1_T(\mathfrak{B}^{\frac dp}_{p,1})}||\chi B^n||_{L^\infty_T(\mathfrak{B}^{\frac dp}_{p,1})}||\nabla \phi||_{\mathfrak{B}^{-\frac dp}_{p',1}}
\\& \quad \ +C||\chi B^n-\chi B||_{L^2_T(\mathfrak{B}^{\frac dp-\frac14}_{p,1})}||\chi u||_{L^2_T(\mathfrak{B}^{\frac dp}_{p,1})}||\nabla \phi||_{\mathfrak{B}^{\frac14-\frac dp}_{p',1}}\rightarrow 0, \quad \mathrm{as} \quad  n\rightarrow \infty,
\end{align*}
which implies $\mathrm{div}(u^n\otimes B^n)$ tends to $\mathrm{div}(u\otimes B)$ in the sense of distributions. Furthermore, we also deduce from Lemma \ref{le2.5} and \eqref{99} that the right-hand terms of the first equation and second equation of \eqref{1.2} belong to $L^1_T(\dot{\mathfrak{B}}^{\frac dp-1}_{p,1})$. Then, by Lemma \ref{le3.2} and \eqref{99}, we have
\begin{align*}
(u,B)\in (\tilde{L}^\infty_T(\dot{\mathfrak{B}}^{\frac dp-1}_{p,1})\cap L^1_T(\dot{\mathfrak{B}}^{\frac dp+1}_{p,1}))^d\times \Big(\tilde{L}^\infty_T(\dot{\mathfrak{B}}^{\frac dp}_{p,1})\Big)^d.
\end{align*}
Finally, following the argument of Theorem 3.19 in \cite{B.C.D}, we can show that $(u,B)\in E^p_T$.\\

\noindent \textbf{Step 4: Uniqueness of the solution. }Assume that $(u^1,B^1)$ and $(u^2,B^2)$ are two solutions of the system with the same initial data. Set $\delta u=u^1-u^2$ and $\delta B=B^1-B^2$. Then, $(\delta u,\delta B)$ satisfies
\begin{equation}\label{4.4}\begin{cases}
\partial_t\delta u-\Delta \delta u=\mathbb{P}\mathrm{div}(-u^1\otimes\delta u-\delta u\otimes u^2+B^1\otimes \delta B+\delta B\otimes B^2),\\
\partial_t\delta B+u^1\cdot\nabla \delta B=\mathrm{div}(u^1\otimes \delta B+\delta u\otimes B^2)-\delta u\cdot\nabla B^2.
\end{cases}\end{equation}
We first claim that $\delta B\in \Big(\tilde{L}^\infty_T(\dot{\mathfrak{B}}^{\frac dp-1}_{p,\infty})\Big)^d$. This claim relies on the following inequality which can be deduced from Lemma \ref{le2.5}, Lemma \ref{le3.1} and $\mathrm{div}(\delta u)=\mathrm{div}(\delta B)=0$:
\begin{align}\label{li}
||\delta B||_{\tilde{L}^\infty_t(\dot{\mathfrak{B}}^{\frac dp-1}_{p,\infty})}&\leq Ce^{CU^1(t)}||\mathrm{div}(u^1\otimes \delta B+\delta u\otimes B^2)-\delta u\cdot\nabla B^2||_{\tilde{L}^1_t(\dot{\mathfrak{B}}^{\frac dp-1}_{p,\infty})}\nonumber\\&\leq Ce^{CU^1(t)}||u^1\otimes\delta B+\delta u\otimes B^2-B^2\otimes\delta u)||_{\tilde{L}^1_t(\dot{\mathfrak{B}}^{\frac dp}_{p,1})}\nonumber\\&\leq Ce^{CU^1(t)}(||\delta u||_{\tilde{L}^1_t(\dot{\mathfrak{B}}^{\frac dp}_{p,1})}||B^2||_{\tilde{L}^\infty_t(\dot{\mathfrak{B}}^{\frac dp}_{p,1})}+||\delta B||_{\tilde{L}^\infty_t(\dot{\mathfrak{B}}^{\frac dp}_{p,1})}||u^1||_{\tilde{L}^1_t(\dot{\mathfrak{B}}^{\frac dp}_{p,1})}),
\end{align}
where $U^1(t)=\int^t_0||u^1||_{\dot{\mathfrak{B}}^{\frac dp+1}_{p,1}}\mathrm{d}\tau$. Using Lemmas $\ref{le2.5}-\ref{le2.6}$ and the fact $\mathrm{div}(\delta u)=\mathrm{div}(\delta B)=0$, we have
\begin{align}\label{4.5}
&\quad \ ||\mathrm{div}(u^1\otimes \delta B+\delta u\otimes B^2)-\delta u\cdot\nabla B^2||_{\tilde{L}^1_t(\dot{\mathfrak{B}}^{\frac dp-1}_{p,\infty})}\nonumber\\&=||\delta B\cdot\nabla u^1+\mathrm{div}(\delta u\otimes B^2-B^2\otimes\delta u)||_{\tilde{L}^1_t(\dot{\mathfrak{B}}^{\frac dp-1}_{p,\infty})}\nonumber\\&\leq C||\delta u||_{\tilde{L}^1_t(\dot{\mathfrak{B}}^{\frac dp}_{p,1})}||B^2||_{\tilde{L}^\infty_t(\dot{\mathfrak{B}}^{\frac dp}_{p,1})}+C||\delta B||_{\tilde{L}^\infty_t(\dot{\mathfrak{B}}^{\frac dp-1}_{p,\infty})}||u^1||_{\tilde{L}^1_t(\dot{\mathfrak{B}}^{\frac dp+1}_{p,1})}.
\end{align}
We apply Lemma \ref{le3.1} and \eqref{4.5} to get for all $t\in[0,T]$,
\begin{align}\label{6.1}
 &\quad \quad||\delta B||_{\tilde{L}^\infty_t(\dot{\mathfrak{B}}^{\frac dp-1}_{p,\infty})}\nonumber\\&\leq Ce^{CU^1(t)}(||\delta u||_{\tilde{L}^1_t(\dot{\mathfrak{B}}^{\frac dp}_{p,1})}||B^2||_{\tilde{L}^\infty_t(\dot{\mathfrak{B}}^{\frac dp}_{p,1})}+||\delta B||_{\tilde{L}^\infty_t(\dot{\mathfrak{B}}^{\frac dp-1}_{p,\infty})}||u^1||_{\tilde{L}^1_t(\dot{\mathfrak{B}}^{\frac dp+1}_{p,1})}).
\end{align}
It follows from Lemma \ref{le2.6} that
\begin{align}\label{4.6}
&\quad \ ||\mathrm{div}(-u^1\otimes \delta u-\delta u\otimes u^2+B^1\otimes \delta B+\delta B\otimes B^2)||_{\tilde{L}^1_t(\dot{\mathfrak{B}}^{\frac dp-2}_{p,\infty})}\nonumber\\&\leq
C||-u^1\otimes\delta u-\delta u\otimes u^2+B^1\otimes \delta B+\delta B\otimes B^2||_{\tilde{L}^1_t(\dot{\mathfrak{B}}^{\frac dp-1}_{p,\infty})}\nonumber\\&
\leq C(||u^1||_{\tilde{L}^2_t(\dot{\mathfrak{B}}^{\frac dp}_{p,1})}+||u^2||_{\tilde{L}^2_t(\dot{\mathfrak{B}}^{\frac dp}_{p,1})})||\delta u||_{\tilde{L}^2_t(\dot{\mathfrak{B}}^{\frac dp-1}_{p,\infty})}\nonumber \\&
\quad +C(||B^1||_{\tilde{L}^1_t(\dot{\mathfrak{B}}^{\frac dp}_{p,1})}+||B^2||_{\tilde{L}^1_t(\dot{\mathfrak{B}}^{\frac dp}_{p,1})})||\delta B||_{\tilde{L}^\infty_t(\dot{\mathfrak{B}}^{\frac dp-1}_{p,\infty})}.
\end{align}
Then, we infer from Lemma \ref{le3.2} and \eqref{4.6} that for any $t\in[0,T]$,
\begin{align}\label{100}
&\quad \ ||\delta u||_{\tilde{L}^1_t(\dot{\mathfrak{B}}^{\frac dp}_{p,\infty})}+||\delta u||_{\tilde{L}^2_t(\dot{\mathfrak{B}}^{\frac dp-1}_{p,\infty})} \nonumber
\\&\leq C||(u^1,u^2)||_{\tilde{L}^2_t(\dot{\mathfrak{B}}^{\frac dp}_{p,1})}||\delta u||_{\tilde{L}^2_t(\dot{\mathfrak{B}}^{\frac dp-1}_{p,\infty})}+C(||B^1||_{\tilde{L}^1_t(\dot{\mathfrak{B}}^{\frac dp}_{p,1})}+||B^2||_{\tilde{L}^1_t(\dot{\mathfrak{B}}^{\frac dp}_{p,1})})||\delta B||_{\tilde{L}^\infty_t(\dot{\mathfrak{B}}^{\frac dp-1}_{p,\infty})}.
\end{align}
We take $\tilde{T}$ small enough such that $||(u^1,u^2)||_{\tilde{L}^2_{\tilde{T}}(\dot{\mathfrak{B}}^{\frac dp}_{p,1})\cap\tilde{L}^1_{\tilde{T}}(\dot{\mathfrak{B}}^{\frac dp+1}_{p,1})}\ll 1/C$. Thus, combining  \eqref{6.1} and \eqref{100}, we infer that for any $t\in[0,\tilde{T}]$,
\begin{align}\label{6.2}
||\delta u||_{\tilde{L}^1_t(\dot{\mathfrak{B}}^{\frac dp}_{p,\infty})}\leq C\int^t_0||\delta B||_{\dot{\mathfrak{B}}^{\frac dp-1}_{p,\infty}}(||B^1||_{\dot{\mathfrak{B}}^{\frac dp}_{p,1}}+||B^2||_{\dot{\mathfrak{B}}^{\frac dp}_{p,1}})\mathrm{d} \tau,
\end{align}
and
\begin{align}\label{6.3}
||\delta B||_{\tilde{L}^\infty_t(\dot{\mathfrak{B}}^{\frac dp-1}_{p,\infty})}\leq Ce^{CU^1(t)}||\delta u||_{\tilde{L}^1_t(\dot{\mathfrak{B}}^{\frac dp}_{p,1})}||B^2||_{\tilde{L}^\infty_t(\dot{\mathfrak{B}}^{\frac dp}_{p,\infty})}.
\end{align}
From Lemma \ref{le2.3}, it follows that
\begin{align*}
||\delta u||_{L^1_t(\dot{\mathfrak{B}}^{\frac dp}_{p,1})}\leq C||\delta u||_{\tilde{L}^1_t(\dot{\mathfrak{B}}^{\frac dp}_{p,\infty})}\log\Big(e+\frac{||\delta u||_{\tilde{L}^1_t(\dot{\mathfrak{B}}^{\frac dp-1}_{p,\infty})}+||\delta u||_{\tilde{L}^1_t(\dot{\mathfrak{B}}^{\frac dp+1}_{p,\infty})}}{||\delta u||_{\tilde{L}^1_t(\dot{\mathfrak{B}}^{\frac dp}_{p,\infty})}}\Big),
\end{align*}
which together with \eqref{6.2} and \eqref{6.3}  yields that for any $t\in[0,\tilde{T}]$,
\begin{align*}
||\delta u||_{\tilde{L}^1_t(\dot{\mathfrak{B}}^{\frac dp}_{p,\infty})}&\leq A_T\int^t_0||\delta u||_{\tilde{L}^1_\tau(\dot{\mathfrak{B}}^{\frac dp}_{p,\infty})}\log(e+\frac{C_T}{||\delta u||_{\tilde{L}^1_\tau(\dot{\mathfrak{B}}^{\frac dp}_{p,\infty})}})\mathrm{d} \tau,
\end{align*}
where
$$A_T=C \exp\{C||u^1||_{L^1_T(\dot{\mathfrak{B}}^{\frac dp+1}_{p,1})}\}||B^2||_{\tilde{L}^\infty_T(\dot{\mathfrak{B}}^{\frac dp}_{p,1})}(||B^1||_{\tilde{L}^\infty_T(\dot{\mathfrak{B}}^{\frac dp}_{p,1})}+||B^2||_{\tilde{L}^\infty_T(\dot{\mathfrak{B}}^{\frac dp}_{p,1})}),$$
and
$$C_T=||\delta u||_{\tilde{L}^1_T(\dot{\mathfrak{B}}^{\frac dp-1}_{p,\infty})}+||\delta u||_{\tilde{L}^1_T(\dot{\mathfrak{B}}^{\frac dp+1}_{p,\infty})}.$$
Note that $A_T$ is integrable on $[0,T]$, and
\begin{align*}
\int^1_0\frac{1}{r\log(e+C_Tr^{-1})}\mathrm{d} r=+\infty.
\end{align*}
An application of Osgood's lemma yields that $(\delta u,\delta B)=(0,0)$ on $[0,\tilde{T}]$, and a continuity argument ensures that $(u^1,B^1)=(u^2,B^2)$ on $[0,T]$. Therefore, combining the above four steps, we complete the proof of Theorem \ref{th1.1}.

\vspace*{1em}
\noindent\textbf{Acknowledgements.}  Li and Yin were
partially supported by NNSFC (No.11671407 and No.11271382),  FDCT (No. 098/2013/A3), Guangdong Special Support Program (No. 8-2015), and the key project of NSF of  Guangdong province (No. 2016A030311004). Tan was partially supported by NNSFC (No. 11301174). The authors thank the referee for valuable comments and suggestions.

\end{document}